\documentclass[letterpaper,conference]{IEEEtran}
\IEEEoverridecommandlockouts
\usepackage{cite}
\usepackage{amsmath,amssymb,amsfonts,amsthm}
\usepackage{algorithmic}
\usepackage{graphicx}
\usepackage{textcomp}
\usepackage{xcolor,url}
\def\BibTeX{{\rm B\kern-.05em{\sc i\kern-.025em b}\kern-.08em
    T\kern-.1667em\lower.7ex\hbox{E}\kern-.125emX}}

\usepackage{pdfsync} 
\usepackage{epsfig}
\usepackage{epstopdf}
\usepackage{amsmath} 
\usepackage{amssymb}
\usepackage{color} 
\usepackage{cite} 
\usepackage{array}
\usepackage{mathrsfs}
\usepackage{cite}

\newcommand{\citep}[1]{\cite{#1}}

\newcommand{\tr}{\mathrm{tr}}

\newcommand{\rank}{\mathrm{rk}}

\newcommand{\field}[1]{\mathbb{#1}} 
\newcommand{\R}{\field{R}}

\newcommand{\Z}{\field{Z}}

\usepackage{nomencl}
\makenomenclature

\usepackage{makeidx}
\makeindex

\usepackage{algorithm}
\usepackage{algorithmic}

\newtheorem{prop}{Proposition}

\newtheorem{ass}{Assumption} 
\graphicspath{{Matlab/}}

\newcommand{\Ree}{\text{Re}}
\newcommand{\Imm}{\text{Im}}

\newcommand{\kEigs}{\mathrm{kEigs}}
\newcommand{\addBlocks}{\mathrm{addBlocks}}
\newcommand{\lift}{\mathrm{lift}}
\newcommand{\blocks}{\mathrm{blocks}}

\newcommand{\get}{\leftarrow}

\usepackage{placeins}

\begin{document}

\title{A Fine-Grained Variant \\ of the Hierarchy of Lasserre \thanks{Contact author:  jakub.marecek@ie.ibm.com;  Wann-Jiun Ma is presently with Stylyze LLC in Sammamish, Washington, United States.}
}
\author{\IEEEauthorblockN{Wann-Jiun Ma}\IEEEauthorblockA{
\textit{IBM Research -- Ireland}\\
Dublin D15, Ireland
} \and \IEEEauthorblockN{Jakub Marecek} \IEEEauthorblockA{
\textit{IBM Research -- Ireland}\\
Dublin D15, Ireland
} \and \IEEEauthorblockN{Martin Mevissen}\IEEEauthorblockA{
\textit{IBM Research -- Ireland}\\
Dublin D15, Ireland
}
}
 
\maketitle

\begin{abstract}
There has been much recent interest in hierarchies of progressively stronger 
convexifications of polynomial optimisation problems (POP). These often 
converge to the global optimum of the POP, asymptotically, 
but prove challenging to solve beyond the first level in the hierarchy for modest instances.
We present a finer-grained variant of the Lasserre hierarchy, together with first-order methods
for solving the convexifications, which allow for efficient warm-starting
with solutions from lower levels in the hierarchy.
\end{abstract}


\section{Introduction}
There has  been much recent interest in efficient solvers for  
polynomial optimisation 
and semidefinite programming (SDP) relaxations therein. 
Much of this interest has been motivated by the work of Lavaei and Low \citep{lavaei2012zero} 
on relaxations of alternating-current optimal power flows (ACOPF), 
an important steady-state problem in power systems engineering \cite{Molzahn2019}. 
For ACOPF, Ghaddar et al.\cite{Ghaddar2015} have shown that the relaxation
of Lavaei and Low is the first level in a number of hierarchies of SDP relaxations, 
including that of Lasserre \cite{Lasserre1},
whose optima converge to the global 
optimum of the polynomial optimisation problem (POP).
The higher levels of the hierarchies present a considerable computational challenge both in ACOPF and other POP.
This is due to the dimension of the relaxation, the super-cubic complexity 
 of traditional interior-point methods for solving SDPs,
as well as their limited warm-start capabilities.
Although a number of hierarchies of second-order cone programming 
problems \cite{7038397,7232429,7285718,7039413} have been studied recently,
for which the interior-point methods are better developed,
the relatively weaker relaxations require yet larger instances to
be solved to provide a strong bound, and have not been proven to dominate
the SDP-based approaches in practice.

In this paper, we suggest that a finer-grained hierarchy of SDP relaxations 
may make sense, if it is accompanied by a method capable of 
``warm-starting,'' i.e., the use of a solution for one relaxation
in speeding up the solution of a stronger relaxation.
Specifically, we:
\begin{itemize}
\item Present a variant of Lassere's hierarchy of SDP relaxations of a POP,
      where localising matrices are added one by one.
\item Design a first-order method for solving these relaxations, which
allows for an efficient warm-starting, and employs a closed-form step in the
coordinate descent.
\item Study the conditions of asymptotic convergence of the approach to the global 
optimum of the POP. 
\end{itemize}
We hope these contributions, alongside
\cite{MarecekTakac2015,LiuMarecekTakac2015,Marecek2016}, 
could spur further interest in first-order methods
for convex relaxations in polynomial optimisation.

\section{Notation and Related Work}

\subsection{Notation}
Let $\R$ and $\Z_+$ be the real numbers and non-negative integers, respectively. Denote by $j$ the imaginary unit of a complex number. The matrix transpose operator and trace operator are denoted by $^T$ and $\text{tr}$, respectively. $\Ree$ and $\Imm$ are the operators that returns the real and imaginary parts of a complex number. Let $e_k$ be the $k$-th standard basis vector in $\mathbb{R}^{N},~N\in \Z_+$. Denote by $\R[x]$ the set of real-valued multivariate polynomials in $x_i,~i=1,...,n$, where $n \in \Z_+$. A polynomial $f \in \R[x]$ is represented as $f(x) := \sum_{\alpha\in \mathcal{F}}c(\alpha)x^{\alpha}$, where $\mathcal{F} \subset \Z_+^n$, and $c(\alpha),~\alpha\in \mathcal{F}$ are the corresponding real coefficients, and $x^{\alpha} := x_1^{\alpha_1}x_2^{\alpha_2}\cdots x_n^{\alpha_n},~\alpha\in \mathcal{F}$ are the corresponding monomials. The support of a polynomial $f\in \R[x]$ is defined by $\text{support}(f):=\{\alpha\in \mathcal{F}\mid c(\alpha)\neq0\}$. The degree of a polynomial $f\in \R[x]$ is defined by $\text{degree}(f):=\text{max}\{\sum_{i=1}^{n}\alpha_i\mid \alpha \in \text{support}(f)\}$. For a non-empty finite set $\mathcal{G} \subset \Z_+^{n}$, $\R[x,\mathcal{G}]:= \{f\in \R[x]  \mid \text{support}(f) \subset \mathcal{G}\}$. $\R[x,\mathcal{G}]^2$ is the set of the SOS polynomials in $\R[x,\mathcal{G}]$. 
  $S(\mathcal{G})$ is the set of $|\mathcal{G}|\times |\mathcal{G}|$ symmetric matrices and $S_+(\mathcal{G})$ is the set of positive semidefinite matrices in $S(\mathcal{G})$ with coordinates $\alpha \in \mathcal{G}$.  $u(x,\mathcal{G})$ is a $|\mathcal{G}|$-dimensional column vector consisting of element $x^\alpha,~\alpha \in \mathcal{G}$. The set $\mathcal{A}_w^C$ is defined as $\mathcal{A}_w^C:=\{\alpha \in \Z_+^n  \mid \alpha_i = 0,~i \notin C, ~\sum_{i\in C}\alpha_i \leq w \}$, for every $C\subset \{1,...n\},~w\in \Z_+$.

\subsection{Polynomial Optimization Problem}\label{s2_ch5}

Let us consider a polynomial optimisation problem:
\begin{equation}\label{pp}
\begin{array}{ll}
\min_{x} & ~ f_0(x) \tag{PP} \\
\mathrm{s.t.}~& ~f_k(x)\geq 0,~k=1,...,m, \notag
\end{array}
\end{equation}
where $x \in \R^n$ is the decision variable, 
the objective and the constraints are defined in terms of are  defined in 
terms of multi-variate polynomials $f_k$, for $k=0,...,m$, in $x \in \R^n.$
A number of approaches have been proposed for solving polynomial optimization problems, including spatial branch-and-bound techniques \cite{belotti2009} and branch-and-reduce \cite{Tawarmalani2005}, cutting plane methods \cite{buchheim2013}, and moment and sum-of-squaress methods \cite{Lasserre1}. Please see \cite{burer2012,Belotti2013,Handbook} for detailed surveys. 
At the same time, noticed that no method can be unconditionally finitely
convergent, as per the solution to Hilbert's tenth problem \cite{Matiyasevich}.

\subsection{Moment-based Methods}

Moment-based methods are a popular approach to solving POP (\ref{pp}), based on the work of Lasserre \citep{Lasserre1}. Given $S\in\R^n$, denote by $\mathcal{P}_w(S)$ the cone of polynomials of degree at most $w$ that are non-negative over $S$. We use $\Sigma_w$ to denote the cone of polynomials of degree at most $d$ that are sum-of-square of polynominals.
Using $\mathcal{G}=\{f_k(x): k=1,\dots,m \}$ and denoting $S_{\mathcal{G}}=\{x \in \mathbb{R}^n :  f(x) \geq 0, \; \forall f \in \mathcal{G}\}$ the basic closed semi-algebraic set defined by $\mathcal{G}$, we can rephrase POP (\ref{pp}) as
\begin{align}\label{pp-d}
\max \quad & \varphi & \mbox{s.t.  }& f(x)-\varphi  \geq 0 \quad \forall \: x \in S_G, \notag\\
= \max \quad & \varphi & \mbox{s.t.  } &f(x)-\varphi \in \mathcal{P}_w(S_{\mathcal{G}}).  
\end{align}
Problem (\ref{pp-d}) is referred to as $[$PP-D$]$.
Although [PP-D] is a conic problem, it is not known how to optimise over the cone $\mathcal{P}_w(S_{\mathcal{G}})$ efficiently. Lasserre \citep{Lasserre1} introduced a hierarchy of SDP relaxations corresponding to liftings of polynomial problems into higher dimensions. In the hierarchy of SDP relaxations, one convexifies the problem, obtains progressively stronger relaxations, but the size of the SDP instances soon becomes computationally challenging. 
Under assumptions slightly stronger than compactness, the optimal values of these problems converge to the global optimal value of the original problem, (\ref{pp}). 

The approximation of $\mathcal{P}_w(S_{\mathcal{G}})$
used by Lasserre \citep{Lasserre1} is the cone $\mathcal{{K}}^w_{\mathcal{G}}$, where
\begin{equation}
\mathcal{{K}}^{w}_{\mathcal{G}} =  \Sigma_{w}+\sum_{k=1}^{m}f_k(x) \Sigma_{w-\deg(f_k)}, 
\end{equation}
and $w \geq d$. 
The corresponding optimization problem over $S$ can be written as:
\begin{equation}\label{eq-Lass}
\begin{split}
\max_{\varphi, \sigma_i(x)} \:  & \varphi  \\
\mbox{s.t. } & f(x)-\varphi= \sigma_0(x)+\sum_{i=1}^{m} \sigma_i(x)f_k(x)\\
&\sigma_0(x) \in \Sigma_{w}, \: \sigma_i(x) \in \Sigma_{w-\deg(f_k)}.  
\end{split}
\end{equation}
Problem (\ref{eq-Lass}) is referred to as $[$PP-H$_w$$]^*$. $[$PP-H$_w$$]^*$ can be reformulated as a semidefinite optimization problem. We denote the dual of [PP-H$_w$]$^*$ by [PP-H$_w$]. The computational cost of the problem clearly depends on both the degree of the polynomials, $w$, and the dimension of the problem. Both the number of constraints and their dimensions can be large when numerous variables of the POP are involved in high-degree polynomial expressions. 

\subsection{Dynamic Generation of Relaxations}

Recently, there has been a considerable effort 
\cite{Kleniati2010,kleniati2010partitioning,kleniati2010decomposition,GhaddarMathProgram,Ghaddar2015,Wittek2015,Campos2017,CAMPOS2018,CAMPOS201932,hall2018optimization} focussed on the design and implementation of decomposition algorithms for solving relaxations [PP-H$_w$] and algorithms for improving those relaxations.
For the hierarchy [PP-H$_w$], Ghaddar et al. \cite{ghaddar2011IPCO,GhaddarMathProgram,Ghaddar2015} have proposed a method, which generated the most violated constraint
in each iteration.
When tested on examples in dimension up to 10 by Ghaddar et al. \cite{Ghaddar2015}, this performs well, but due to the limitations of current SDP solvers, does not scale much further.
A similar approach is being considered by Molzahn and Hiskens \citep{MolzahnHiskens2015a} in a power-system specific heuristic,
and by Chen \citep{chen2012maximum} for general polynomial optimisation problems.
Kleniati et al. \cite{Kleniati2010,kleniati2010partitioning,kleniati2010decomposition} study a variant of Bender's decomposition,
where the moment variable is decomposed into blocks, such that the constraints particular to each block are considered in isolation
in one sub-problem, and the constraints spanning multiple blocks are considered in a so called master problem.
In \cite{kleniati2010partitioning} convergence to the optimum of the polynomial optimisation problem is studied,
with details of run-time provided on examples in dimension up to 10.
Wittek \cite{Wittek2015} describes a mixed-level relaxation, where monomials can be added arbitrarily, but without an algorithmic approach for their addition.
In the context of relaxations of partial differential equations (PDEs), the so-called prolongation operators are used routinely.
Campos \cite{Campos2017,CAMPOS2018,CAMPOS201932} has translated this work to SDP relaxations of PDEs and beyond.
Hall \cite{hall2018optimization} describes a ``Sum of Squares Basis Pursuit'' using linear or second-order cone programming, but
also shows \cite[Proposition 3.5.5]{hall2018optimization} the approach is not convergent.
Our hope is to improve upon this state of the art.

\FloatBarrier

\section{The Hierarchy}\label{s4_ch5}


Let us consider the dual of Problem (\ref{eq-Lass}), i.e., the semidefinite programming relaxation obtained by Lasserre in the method of moments. Following Chapter 6 in \cite{lasserre2015introduction}, we can write it as: 
\begin{align}\label{sdpa}
\underset{y}{\inf }& ~F(y) \tag{PP-H$_w$} & \\
\mathrm{s.t. }& ~ M_w(y) \succeq 0 \label{moment} & \\
              & ~ M_{w-v_j}(f_k y) \succeq 0 & \forall \; k = 1, \ldots, m \label{localising} \\
              & ~ y_0 = 1
\end{align}  
where $F$ is a linear functional, $M_w(y)$ is called the moment matrix and 
$M_w(y) = \sum{\alpha} y_\alpha C_\alpha^0$ for some appropriate real matrix $C_\alpha^0$. 
$M_{w-v_j}(f_k y)$ is called a localising matrix and 
$M_{w-v_j}(f_k y) = \sum{\alpha} y_\alpha C_\alpha^0$ 
some appropriate real matrix $C_\alpha^0$ for each inequality $k = 1, \ldots, m$.

It is clear we could replace constraints \eqref{moment} and \eqref{localising} with 
a constraint on a single block-diagonal matrix to be positive definite,
where the blocks on the diagonal would be $M_w(y)$ and $M_{w-v_j}(f_k y), k = 1, \ldots, m$.
Let us denote these blocks $B_w$ and let us use the notation 
$\prod_{b \in \mathcal B_w} b$  for the formation of the block-diagonal matrix, with the blocks
taken in arbitrary order.
We can use a sub-set of blocks $B_q \subseteq \mathcal B_w$ to be considered in iteration $q$. Let us also consider the complement $\bar B_q$ of the block $B_q$, i.e., $B_q \cup \bar B_q = \mathcal B_w$ for the $w$ current at $j$.
Our hierarchy is simply based on the SDP relaxation $[\textrm{PP-H}_w]$ parametrised by the choice $B_q$ of the blocks at the $q$-th iteration within the $w$-th level of relaxation: 
\begin{align}\label{sdpr}
\underset{y}{\inf }& ~L(y) \tag{R($w,B_q$)} \\
\mathrm{s.t. }& ~ \prod_{b \in \mathcal B_q} b \succeq 0 \\
              & ~ y_0 = 1.
\end{align}  

It is not clear what blocks $B_q \subseteq \mathcal B_w$ to consider, though. In determining those, we use two maps:
\begin{itemize}
\item $F: \mathcal B \rightarrow [m]$ maps blocks to constraints $f_j$ of the polynomial optimization problem
\item $G: [m] \rightarrow N$ maps the constraints to the variables $i\in \mathcal{N}$ in the polynomial optimization problem. 
\end{itemize}
The composite mapping $G \circ F : \mathcal B \rightarrow N$ hence maps the blocks to the variables of the polynomial optimization problem and describes the relationship between the POP and the SDP. 
Notice that $F, G$ are known.
One can construct $F, G$ in the process of formulating the moment and localising
matrices.
Alternatively, one can use the simplistic procedure for obtaining $F, G$, such as Algorithm \ref{blocks}.

\section{An Algorithm}\label{ch5_sec_algo}

Let us describe the complete algorithm for solving the polynomial optimisation problem, based on:
\begin{itemize}
\item The moment-based relaxations (\ref{sdpr}) suggested above.
\item A novel ``all violated'' block-addition rule, for picking suitable blocks  to add to the relaxation, building upon  the ``most violated'' \citep{chen2012maximum} and ``power mismatch'' \citep{MolzahnHiskens2015a} block-addition rules. 
\item The augmented Lagrangian approach, an optimisation strategy studied since the 1950s, e.g., by Hestenes \citep{Hestenes1969} and Powell  \citep{Powell1969}. 
\item The block-wise additivity of the augmented Lagrangian in the moment-based relaxations, as outlined above and developed further below.
\item A novel closed-form step for the augmented Lagrangian approach, as applied to the  moment-based relaxations, as explained below.
\end{itemize}
The ingredients, which are truly novel, are the block-separable augmented Lagrangian and the closed-form step,
but we argue that the overall algorithm design should also be of interest.
Subsequently, we prove the convergence of the algorithm to the global optimum in the following section.

\subsection{The Overall Algorithm}\label{s5_ch5}


\begin{algorithm}[t!]
\caption{Optimization with the Fine-Grained Variant of the Hierarchy of Lasserre} 
\label{bcd_sdp_sh}
\begin{algorithmic}[1]
\item [Input: (PP), i.e. objective $f_0$ and constraints $f_j,~j=1,...,m$]
\STATE $X^0\in \tilde{S}_+$, $Z^0\in \tilde{S}_+$, $y^0\in \R^{|\tilde{\mathcal{F}}|}$
\STATE $D \get \max_i \deg(f_i)$, $d \get \lceil D/2 \rceil$
\STATE $B_0 \get$ a list containing all the blocks in (PP-H$_d$) \label{ini_B}
\STATE $w \get d + 1$\label{ini_w}
\STATE $q\get 0$ 
\REPEAT \label{outer_1}
\label{l:rankcriteria1}
\STATE $(\bar{B_q},F,G) \get \blocks(w)$, i.e. the list of all blocks at level $w$ \label{l:barB}
\REPEAT
\label{l:rankcriteria2}
\STATE $k\get 0$
\color{black}\WHILE{\text{Cauchy criteria are not satisfied}} \label{l:BCD}
    \STATE $l\get 0$, $y^{kl}\get y^k$
\WHILE{\text{Cauchy criteria are not satisfied}}
\FOR{coordinate $\alpha= 1,2,...$ \textbf{ in parallel } } 
\STATE compute $a_1(\alpha),~a_2(\alpha)$ using \eqref{aa1}, \eqref{aa2}
\STATE $y_i^{kl}\get  \frac{-a_2(\alpha)}{2a_1(\alpha)}$
\ENDFOR
\STATE $l\get l+1$
 \ENDWHILE
\STATE $y^{k+1}\get y^{kl}$
\FOR{block $b \in B_q$ \textbf{ in parallel }}
\STATE $V_b^{k+1} \get \bar{M}_b(y^{k+1})+ \mu Z_b^{k}$
\STATE $(E^r_{b+}, V^r_{b+}) \get \kEigs(V_b^{k+1})$ \label{l:eigs}, 
\STATE $X_b^{k+1}\get V^r_{b+}E^r_{b+}V^{rT}_{b+}$
\STATE $Z_b^{k+1}\get V_b^{k+1}-X_b^{k+1}$
\ENDFOR 
\STATE $k\get k+1$    
\ENDWHILE
\label{l:BCDEnd}

\STATE $(B_{q+1}, \bar B_{q+1}, s) \get \addBlocks(X^{k}, B_q, \bar B_q,F,G)$ \label{l:cuttingsurfaces}
\STATE 
       $q\get q+1$
\IF{$s = 0$ \textbf{ and } $\bar B_{q} \not = \emptyset$ \textbf{and} \textrm{approximateObj}$(f_0, \mathcal{F}, \tilde{c}, w, X^k, y^k)$} \label{convergenceAchieved}
  \STATE \textbf{return} $(X^k, y^k, Z^k)$
\ENDIF
\UNTIL{$s = 0$ \textbf{ or } $\bar B_{q} = \emptyset$ }
\color{black}\STATE $w \get w + 1$
\color{black}\STATE $(X^{k}, Z^{k}, y^{k}) \get \lift(X^{k}, Z^{k}, y^{k})$\label{lift}\label{outer_2} 
\UNTIL{
  \textrm{flatExtension}$(D, d, w, X^k, y^k)$ 
}
\STATE \textbf{return} $(X^k, y^k, Z^k)$
\end{algorithmic}
\end{algorithm} 

Algorithm \ref{bcd_sdp_sh} captures the key algorithm schema, with details elaborated in Algorithms~\ref{blocks}--\ref{approximate}.
In Algorithm \ref{bcd_sdp_sh}, we first initialise $d$ to one half of the maximum degree involved in (PP), rounded up, and 
 $B$ is set to a list of blocks considered in (PP-H$_d$).
The main loop of the algorithm has $w$ as the counter, which denotes the order of the relaxation, from which we are adding blocks in that particular iteration.
On line~\ref{l:barB}, $\mathrm{blocks}(w)$ constructs the mapping between the buses and the blocks in the $w$-th level of relaxation. 
$\bar B$ is a list of blocks from level $w$ of the hierarchy, which has not been part of the SDP relaxation yet.
$F(j)$ is a mapping between the $j$-th polynomial constraint and the corresponding block. $G(j)$ is a mapping between the $j$-th constraint and the  corresponding buses. 
On lines~\ref{l:BCD}--\ref{l:BCDEnd}, we solve the SDP relaxation using the block-coordinate descent method, where the blocks correspond to the blocks from (PP-H$_{w-1}$) and $B$. 
On line~\ref{l:eigs}, $\kEigs(X)$ is a procedure, which obtains the largest $r$ eigenvalues with the associated eigenvectors of $X$. This could be simply implemented using a spectral decomposition of $X$, where all but the largest $r$ eigenvalues and the associated eigenvectors are discarded. 
Alternatively, one could use the power method repeatedly with deflation, or the extensions of the power method extracting multiple eigenvalues in one pass. On line~\ref{l:cuttingsurfaces}, $\addBlocks$ is a procedure, which computes the projection of $X$ onto rank one matrices, i.e. the closest rank one matrix to $X$ with respect to the Frobenius norm, and one-by-one verifies whether the constraints are satisfied up to the accuracy $\epsilon = 10^{-5}$. Whenever a constraint at bus $v$ is not satisfied, the blocks in $\bar B$ corresponding to bus $v$ are added. The number of blocks
added to $B$ in this run of $\addBlocks$ is output into $s$.
The Cauchy convergence criteria are set to $10^{-3}$, on both lines~\ref{l:rankcriteria1} and \ref{l:rankcriteria2}.
Finally, at the end of the outer loop on line \ref{lift}, $\lift(X, Z, y)$ is a procedure, which lifts the two matrices $X, Z$ and vector $y$ to a higher dimension.

\subsection{The Choice of Blocks}

Next, we suggest the procedure for the addition of blocks with constraints violated at the current relaxation to the relaxation. 
Alternatively, one can see that as the removal of the redundant blocks in the relaxation.

\begin{algorithm}[t!]
\caption{blocks}\label{blocks}
\begin{algorithmic}[1]
\item [Input: level $w$ in the hierarchy (PP-H$_w$)]
\STATE $\bar{B} \get $ all the blocks in relaxation (PP-H$_w$)\label{allBlocks},  
\FOR{$j=1,...,m$}
\STATE$\textrm{perturb coefficients of the constraint}~~f_{j}$\label{change_constraint}\label{F_1}
\FOR{$\alpha=1,...,|\mathcal{F}_w|$}
\STATE $\tilde{y}_{\alpha} \get \mathrm{rand}(1,1)$
\ENDFOR
\STATE $M(\tilde{y})^{-} \get M(0)+\sum_{\alpha\in \tilde{\mathcal{F}_{w}}}M(\alpha)\tilde{y}_{\alpha}
- (\tilde{M}(0)+\sum_{\alpha\in \tilde{\mathcal{F}_w}}\tilde{M}(\alpha)\tilde{y}_{\alpha})$
\STATE $F(j)\get$ map the remaining elements in $M(\tilde{y})^{-}$ to the corresponding block$\label{F_2}$
\STATE $G(j) \get$ all the buses involved in the constraint~~$f_j$\label{G}
\ENDFOR
\STATE \textbf{return} $(\bar{B}, F,G)$
\end{algorithmic}
\end{algorithm} 


In Algorithm \ref{blocks}, on line \ref{allBlocks}, $\bar{B}$ is a list which contains all the blocks in the $w$-th level of relaxation. On line \ref{F_1} to \ref{F_2}, we construct the mapping $F(j)$, for every $j=1,...,m$. On line \ref{G}, we construct the mapping $G$.

\begin{algorithm}[t!]
\caption{addBlocks}\label{addBlocks}
\begin{algorithmic}[1]
\item [Input: $X,~B,~\bar{B},~F,~G$]
\STATE $s \get 0$, $\epsilon \get 10^{-5}$, $I \get \emptyset$
\STATE $M_2 \get$ submatrix of $X$ corresponding to the second-order monomials 
\STATE $(E, V) \get \kEigs(M_2, 1)$\label{kEigsW} 
\STATE $x \get \textrm{chol}(VEV^T)$\label{Wrank1}
\FOR{constraint $k=1,...,m$} \label{mismatch_diff_1}
\IF{$|\min(f_k(x), 0)| > \epsilon $}
\STATE $I \get I+ \{k\}$
\ENDIF
\ENDFOR \label{mismatch_diff_2}
\FOR{$j = 1,...,m$}\label{remove_1}
\IF{$\{G(j)\bigcap I\} \neq \emptyset$}
\STATE $B \get B \cup \{F(j)\}$
\STATE $s \get s+1$\label{supdate}
\ENDIF
\ENDFOR\label{remove_2}
\STATE $\bar{B} \get \bar{B}\setminus B$\label{Bbar_update}
\STATE \textbf{return} $(B,\bar{B},s)$
\end{algorithmic}
\end{algorithm} 


In Algorithm \ref{addBlocks}, 
we pick blocks to add to the current relaxation, which correspond to constraints violated by more than $\epsilon$ by the current iterate. 
In order to obtain vector $x$ of (PP) from the current iterate, on line \ref{kEigsW}, $\kEigs$ computes the largest eigenvalue and the associated eigenvector of the matrix corresponding to second-order monomials. Subsequently, we obtain a vector $x$ by Cholesky decomposition of the rank-1 projection on line \ref{Wrank1}. On line \ref{mismatch_diff_1}--\ref{mismatch_diff_2}, we verify if the power constraints are satisfied up to the accuracy $\epsilon= 10^{-5}$. 
On line \ref{remove_1} to \ref{remove_2}, we remove the redundant blocks involved in the violated constraints. 
On line \ref{supdate}, we update the total number of blocks added and removed. On line \ref{Bbar_update}, the blocks added to the list $B$ are removed from the list $\bar{B}$.

\begin{algorithm}[t!]
\caption{flatExtension}\label{flatextension}
\begin{algorithmic}[1]
\item [Input: constants $D, d$, primal-dual pair $X,y$ level $w \in \R$]
  \STATE $M_2, M_3, \ldots, M_w \gets $ sub-matrices of $X$ with moment matrices of the order
  \IF {$\left( \right. (D \le w \textbf{ and } \rank(M_w) = \rank(M_{w - 1})$ \\ 
  \quad \quad \textbf{or} ($d \le w$ \textbf{ and } $\rank(M_w) = \rank(M_{w - d}) \left. \right)$ }
  \STATE \textbf{ return True } 
  \ENDIF
  \STATE \textbf{ return False } 
\end{algorithmic}
\end{algorithm} 

\begin{algorithm}[t!]
\caption{approximateObj}\label{approximate}
\begin{algorithmic}[1]
\item [Input: instance as $f_0, \mathcal{F}, \tilde{c}$, primal-dual pair $X,y$ level $w \in \R$]
  \STATE $M_2, M_3, \ldots, M_w \gets $ sub-matrices of $X$ with moment matrices of the order
  \STATE $z_{SDP} \gets \sum_{\alpha \in \tilde{\mathcal{F}}}\tilde{c}_0(\alpha)y_{\alpha}$
  \STATE $(E, V) \get \kEigs(M_2, 1)$ 
  \STATE $x \get \textrm{chol}(VEV^T)$
  \STATE $z_{POP} \gets f_0(x)$ 
  \IF{ $|z_{SDP} - z_{POP}| \le \epsilon $ }
  \STATE \textbf{ return True } 
  \ENDIF
  \STATE \textbf{ return False } 
\end{algorithmic}
\end{algorithm} 

In Algorithms \ref{flatextension} and \ref{approximate}, we suggest the usual flat extension 
test of convergence, and a comparison of objective-function
values of the current iterate and the vector obtained by Cholesky decomposition of the 
closest rank-1 projection of the submatrix of $X$ corresponding to second-order monomials.
The comparison of objective-function values needs to be considered alongside the satisfaction
of all constraints of (PP), as suggested on Line~\ref{convergenceAchieved} of the main algorithm.

\subsection{The Augmented Lagrangian}

The augmented Lagrangian $L_{\mu}$ of (\ref{sdpa}) is defined by
\begin{equation}\label{lagrangian}
L_{\mu}(Z,y,X) = c^Ty+<Z,\bar{M}(y)-X>+\frac{1}{2\mu}||\bar{M}(y)-X||_F^2,
\end{equation}
where $Z \in \tilde{S}$ is the dual variable and the parameter $\mu$ is a positive real number.  
Typically, this is used in conjunction with the alternating direction method of multipliers (ADMM), which in iteration $k$ computes the updates:
\begin{equation}\label{yupdate}
y^{k+1}=\underset{y}{\min}~L_{\mu}(Z^k,y,X^k),
\end{equation}
\begin{equation}\label{xupdate}
X^{k+1}=\underset{X\succeq 0}{\min}~L_{\mu}(Z^k,y^{k+1},X),
\end{equation}
\begin{equation}\label{zupdate}
Z^{k+1}=Z^k+\frac{\bar{M}(y^{k+1})-X^{k+1}}{\mu}.
\end{equation}  
Notice that we effectively perform a two-block decomposition of the augmented
Lagrangian, rather than the multi-block decomposition, which is known \cite{chen2014direct} to be divergent in some cases, esp. when 
\cite{hong2012linear} the functions involved are not strongly convex.

In our case, we have the following special structure:
\begin{prop}
\label{additive}
The augmented Lagrangian is additive with respect to the blocks.
\end{prop}

\emph{Proof:}
The constraint $X \in \tilde{\mathcal{S}}_+$ is sufficient to ensure that $\bar{M}(y)-X$ has the same block structure of $\bar{M}(y)$. To complete the proof, we  use the fact that the trace of a block diagonal matrix with square blocks is equal to the sum of the traces of the blocks.  

Given the block diagonal structure of $\bar{M}$(y), we are able to decompose the computation block-wise.

The first-order optimality condition for (\ref{xupdate}) yields 
\begin{equation}
-Z^k + \frac{1}{\mu}(X^{k+1}-\bar{M}(y^{k+1}))  = 0,
\end{equation}
and
\begin{equation}\label{v}
X^{k+1} = \bar{M}(y^{k+1})+ \mu Z^{k},~X^{k+1}\succeq 0.
\end{equation}

In order to find the solution of (\ref{v}), spectral decomposition is performed on the matrix
\begin{equation}
V^{k+1} := \bar{M}(y^{k+1})+ \mu Z^{k},
\end{equation} 
and the result is used to formulate $V_+E_+V_+^T$, where $E_+$ contains the nonnegative eigenvalues of the matrix $V^{k+1}$, and the columns of $V_+$ are the corresponding eigenvectors. 

Substitute $\bar{M}(y^{k+1})+ \mu Z^{k}$ by $V^{k+1}$ into (\ref{zupdate}), which yields 
\begin{equation}\label{zupdate2}
Z^{k+1}=V^{k+1}-X^{k+1}.
\end{equation} 

\subsection{Block-Coordinate Descent Method}

Further, we need a method to compute the update \eqref{yupdate} efficiently. Considering \ref{additive}, we suggest:

\begin{prop}
For every $i\in \tilde{\mathcal{F}}$, the first-order optimality conditions for the $i$-th coordinate of $y$ in (\ref{yupdate}) yield:
\begin{equation}\label{y_closed_form}
y_i^{k+1} = \frac{-a_2(i)}{2a_1(i)},
\end{equation}
where
\begin{align}
a_1(i) &= \frac{1}{2\mu}\tr\big{(}M_iM_i^T\big{)},\label{aa1}\\
a_2(i) &= c(i)+\tr(Z^kM_i^T\big{)}+\sum_{j\neq i}\frac{1}{\mu}\tr\big{(}M_jy_jM_i^T\big{)}\notag\\
&-\frac{1}{\mu}\tr\big{(}M(0)M_i^T\big{)}-\frac{1}{\mu}\tr\big{(}M_i^TX^k\big{)}\label{aa2}.
\end{align}
\end{prop}

\emph{Proof:}
Note that $L_{\mu}(Z^k,y_{i},X^k)$ is a quadratic function of $y_{i}$ and algebraic manipulations lead to the result. 

When we update the $\alpha$-th coordinate of $y$, the only blocks of $\bar{M}(y)$ that are required are the ones containing $y_{\alpha}$. Each block of $\bar{M}(y)$ only contains a portion of the variable $y$, which results in speeding up the implementation.

\FloatBarrier

\section{An Analysis}

Algorithm \ref{bcd_sdp_sh} has been designed so as to be convergent. 
In particular:
\begin{itemize} 
\item in the outer-most loop (Lines 3--23), for each $w$, we compute the
optimum of a relxation [OP$_2$-H$_w$]$^*$. 
In the limit of $w$, relaxations [OP$_2$-H$_w$]$^*$ converge
 to the
global optimum of [OP$_2$],
as shown by Ghaddar et al.\ \citep{Ghaddar2015}. 
However, the relaxations are \emph{not} formed explicitly,
but dynamically, by considering the violated blocks.
\item in the block-addition loop (Lines 6--21), specifically, we include
the violated blocks. The finiteness of the loop is given by 
the finiteness of the number of blocks in [OP$_2$-H$_w$]$^*$ for any
finite $w$.
\item in the inner-most loop (Lines 8--19), the dynamically constructed 
relaxation of [OP$_2$-H$_w$]$^*$ is solved by a first-order method.
The convergence is based on a rich history of work on 
the convergence of first-order methods for semidefinite programming.
\end{itemize}

In the analysis, we start from the inner-most loop and proceed outwards,
while formulating the assumptions before we use them. Throughout, we 
use the notion of the dual problem of (\ref{sdpa}), which is
\begin{equation}\label{primal}
\begin{array}{ll}
\underset{Z}{\max}~~-\langle M(0), Z\rangle\\
\mathrm{s.t.}~~~\langle M(i), Z\rangle = c_i,~i=1,...m\\
~~~~~~~~~~~~~~~~Z \succeq 0.
\end{array}
\end{equation}  
 
\begin{ass} \label{Slater_P}
For problem (\ref{primal}), 
the set of feasible solutions is compact and 
there exists $Z_0 \in S^n_{+}$ such that
\begin{eqnarray} \label{eq-slater-p}
&\langle M(i), Z_0\rangle = c_i,~i=1,...m,\quad Z_0 \succ 0.
\end{eqnarray}
\end{ass}

\begin{ass} \label{Slater_D}
For problem (\ref{sdpa}),
the set of feasible solutions is compact and  
there exists $(y_0,X_0) \in \R^m \times \mathcal{S}^n_{+}$
such that 
\begin{eqnarray} \label{eq-slater-d}
\bar{M}(y_0) = X_0 , \; X_0 \succ 0,
\end{eqnarray}
\end{ass}

Notice that if there exist feasible but non-strict solutions of (\ref{eq-slater-p}) or (\ref{eq-slater-d}), by eliminating implicit equality constraints, any feasible but non-strictly feasible
instance of SDP can be reduced to an equivalent strictly feasible instance of SDP, which satisfies Assumption \ref{Slater_P} and \ref{Slater_D},
as per Section 2.5 in \citep{bgfb} and \citep{vb}.

%

Further, with respect to the convergence of the first-order method, in order to simplify the proof, let us assume that:

\begin{ass} 
\label{SmalesModel}
There exists a unique optimum of (\ref{yupdate}) and (\ref{zupdate}),
 which can be computed exactly, just as all other computations 
 are performed exactly.
\end{ass} 

Such assumptions are standard \cite{razaviyayn2013}.
The accuracy benefits from the use of a closed-form formula \eqref{y_closed_form},
considering one may obtain an irrational number even when all coefficients are integral, 
there are rounding errors and the errors may propagate.
We do not perform an analysis of error propagation throughout the computation, 
which renders our analysis closer to the real computation model of Blum et al. \cite{blum2012complexity}, than 
the Turing model customarily used in Computer Science.
This assumption can be relaxed \cite{tappenden2013inexact}, albeit at the
price of a substantially more complex analysis. 

Specifically:


\begin{prop}
\label{convergence-inner-most}
If Assumptions \ref{Slater_P}--\ref{SmalesModel} hold, the sequence $\{Z^k\}$ generated in the inner loop (Lines 8--19) converges
to $Z^*$, where $Z^*$ is an optimal solution of
\ref{primal}, and the augmented Lagrangian
$L_{\mu}(Z^k,y^k,X^k)$ converges to $p^*$ as $k$ goes to infinity, where $p^*=$ \eqref{sdpa} = \eqref{primal}.
\end{prop}

\emph{Proof:}
The global convergence follows from Theorem 8 of \citep{eb},
but could be derived also directly from Theorem 1 of \citep{R76a}
and Theorem 4 of \citep{R76b}.
For the stopping rule on Line 11 of Algorithm  \ref{bcd_sdp_sh},
one could envision each subsequent iteration of the loop on Lines 9--20 using an increased precision, but one can equally well prove convergence with limited precision.

\begin{prop}
\label{convergence-blocks}
If Assumptions \ref{Slater_P}--\ref{SmalesModel} hold, for each $w \ge 3$, the sequence 
of optima of the semidefinite-programming relaxations
\ref{sdpr} generated in block-addition 
loop (Lines 6--21) converges
to \ref{sdpa} 
as $q$ goes to infinity. 
\end{prop}

\emph{Proof:}
The proof follows from Theorem 5 of \cite{KojimaMuramatsu2009};
The convergence is finite for any finite $w$, because the number of blocks is finite \cite{KojimaMuramatsu2009}, for any finite $w$.
Notice we have replaced the assumptions, in line with \cite{Ghaddar2015}.

Overall, we have the convergence result as follows:


\begin{prop}
\label{thm:convergent}
If Assumptions \ref{Slater_P}--\ref{SmalesModel} hold, 
$\inf \{L^w_{\mu}(Z^k,y^k,X^k)\}\rightarrow \min \text{([PP])}$ as $w\rightarrow\infty,~k\rightarrow\infty.$ 
\end{prop}

\emph{Proof:}
The proof combines Propositions~\ref{convergence-inner-most} and 
\ref{convergence-blocks} above. 
The complication is in the stopping rules for the coordinate-descent (Line 11) and  overall solver of the SDP (Line 9), such that we can truly obtain the (possibly irrational) solution of the correct SDP, while adding blocks. 
For the purposes of the proof, we suggest a dynamic stopping rule for Line 9 of Algorithm  \ref{bcd_sdp_sh}, where in each subsequent SDP solved, i.e. each subsequent iteration of the loop on Lines 9--20, the precision required is doubled.
A dynamic stopping rule for Line 11 of Algorithm  \ref{bcd_sdp_sh} can be similarly doubled in each iteration.
This would clearly lead to an arbitrary precision, eventually, and would be able to produce the possibly irrational numbers, asymptotically. 


Notice that for many problems, including the ACOPF, Assumptions are \ref{Slater_P}, \ref{Slater_D} are satisfied for all realistic choices
 of parameters, as the feasible region is compact.
Although the assumptions may be relaxed slightly, 
the work of Matiyasevich \cite{Matiyasevich} suggests that unconditional
finite convergence is impossible.


\section{Conlusions}

The development of practical and globally convergent solvers for polynomial optimisation problems is a major challenge within mathematical optimisation. 
In turn, this poses challenges in convex optimisation, including first-order methods for semidefinite-programming,
and numerical linear algebra, such as the incremental  update of the (truncated) singular value decomposition \cite{Bunch1978,gu1994,Brand2006}.
Overall, we have made preliminary steps towards the development of a method, which is both convergent and efficient in practice, although there are still very distinct limitations. 

\clearpage
\bibliographystyle{abbrv} 
\bibliography{acopf,literature,references}

\begin{thebibliography}{10}

\bibitem{Handbook}
M.~F. Anjos and J.-B. Lasserre, editors.
\newblock {\em Handbook on semidefinite, conic and polynomial optimization},
  volume 166 of {\em International series in operations research \& management
  science}.
\newblock Springer, New York, 2012.

\bibitem{Belotti2013}
P.~Belotti, C.~Kirches, S.~Leyffer, J.~Linderoth, J.~Luedtke, and A.~Mahajan.
\newblock {M}ixed-{I}nteger {N}onlinear {O}ptimization.
\newblock In A.~Iserles, editor, {\em {A}cta {N}umerica}, volume~22, pages
  1--131. Cambridge University Press, 2013.

\bibitem{belotti2009}
P.~Belotti, J.~Lee, L.~Liberti, F.~Margot, and A.~W{\"a}chter.
\newblock Branching and bounds tightening techniques for non-convex minlp.
\newblock {\em Optimization Methods \& Software}, 24(4-5):597--634, 2009.

\bibitem{blum2012complexity}
L.~Blum, F.~Cucker, M.~Shub, and S.~Smale.
\newblock {\em Complexity and real computation}.
\newblock Springer Science \& Business Media, 2012.

\bibitem{bgfb}
S.~P. Boyd, L.~El~Ghaoui, E.~Feron, and V.~Balakrishnan.
\newblock {\em Linear matrix inequalities in system and control theory},
  volume~15.
\newblock SIAM, 1994.

\bibitem{Brand2006}
M.~Brand.
\newblock Fast low-rank modifications of the thin singular value decomposition.
\newblock {\em Linear Algebra and its Applications}, 415(1):20 -- 30, 2006.
\newblock Special Issue on Large Scale Linear and Nonlinear Eigenvalue
  Problems.

\bibitem{buchheim2013}
C.~Buchheim and A.~Wiegele.
\newblock Semidefinite relaxations for non-convex quadratic mixed-integer
  programming.
\newblock {\em Mathematical Programming}, 141(1-2):435--452, 2013.

\bibitem{Bunch1978}
J.~R. Bunch and C.~P. Nielsen.
\newblock Updating the singular value decomposition.
\newblock {\em Numerische Mathematik}, 31(2):111--129, 1978.

\bibitem{burer2012}
S.~Burer and A.~N. Letchford.
\newblock Non-convex mixed-integer nonlinear programming: a survey.
\newblock {\em Surveys in Operations Research and Management Science},
  17(2):97--106, 2012.

\bibitem{Campos2017}
J.~Campos and P.~Parpas.
\newblock A multigrid approach to sdp relaxations of sparse polynomial
  optimization problems.
\newblock {\em SIAM Journal on Optimization}, 28(1):1--29, 2018.

\bibitem{CAMPOS2018}
J.~S. Campos.
\newblock {\em A multigrid approach to SDP relaxations of sparse polynomial
  optimization problems}.
\newblock PhD thesis, Imperial College London, 2018.

\bibitem{CAMPOS201932}
J.~S. Campos, R.~Misener, and P.~Parpas.
\newblock A multilevel analysis of the lasserre hierarchy.
\newblock {\em European Journal of Operational Research}, 277(1):32--41, 2019.

\bibitem{chen2012maximum}
B.~Chen, S.~He, Z.~Li, and S.~Zhang.
\newblock Maximum block improvement and polynomial optimization.
\newblock {\em SIAM Journal on Optimization}, 22(1):87--107, 2012.

\bibitem{chen2014direct}
C.~Chen, B.~He, Y.~Ye, and X.~Yuan.
\newblock The direct extension of admm for multi-block convex minimization
  problems is not necessarily convergent.
\newblock {\em Mathematical Programming}, pages 1--23, 2014.

\bibitem{eb}
J.~Eckstein and D.~P. Bertsekas.
\newblock On the {D}ouglas--{R}achford splitting method and the proximal point
  algorithm for maximal monotone operators.
\newblock {\em Math. Program.}, 55(1-3):293--318, 1992.

\bibitem{Ghaddar2015}
B.~Ghaddar, J.~Mare{\v c}ek, and M.~Mevissen.
\newblock Optimal power flow as a polynomial optimization problem.
\newblock {\em IEEE Trans. Power Syst.}, 31(1):539--546, 2016.

\bibitem{GhaddarMathProgram}
B.~Ghaddar, J.~Vera, and M.~Anjos.
\newblock A dynamic inequality generation scheme for polynomial programming.
\newblock {\em Mathematical Programming}, pages 1--37, 2015.

\bibitem{ghaddar2011IPCO}
B.~Ghaddar, J.~C. Vera, and M.~F. Anjos.
\newblock An iterative scheme for valid polynomial inequality generation in
  binary polynomial programming.
\newblock In {\em Integer Programming and Combinatoral Optimization}, pages
  207--222. Springer Berlin Heidelberg, 2011.

\bibitem{gu1994}
M.~Gu and S.~C. Eisenstat.
\newblock A stable and fast algorithm for updating the singular value
  decomposition, 1994.
\newblock Yale University Technical Report.

\bibitem{hall2018optimization}
G.~Hall.
\newblock {\em Optimization over Nonnegative and Convex Polynomials With and
  Without Semidefinite Programming}.
\newblock PhD thesis, Massachusetts Institute of Technology, 2018.
\newblock arXiv preprint arXiv:1806.06996.

\bibitem{Hestenes1969}
M.~R. Hestenes.
\newblock Multiplier and gradient methods.
\newblock {\em Journal of optimization theory and applications}, 4(5):303--320,
  1969.

\bibitem{hong2012linear}
M.~Hong and Z.-Q. Luo.
\newblock On the linear convergence of the alternating direction method of
  multipliers.
\newblock {\em arXiv preprint arXiv:1208.3922}, 2012.

\bibitem{kleniati2010decomposition}
P.~Kleniati, P.~Parpas, and B.~Rustem.
\newblock Decomposition-based method for sparse semidefinite relaxations of
  polynomial optimization problems.
\newblock {\em Journal of optimization theory and applications},
  145(2):289--310, 2010.

\bibitem{Kleniati2010}
P.~M. Kleniati.
\newblock {\em Decomposition schemes for polynomial optimisation, semidefinite
  programming and applications to nonconvex portfolio decisions}.
\newblock PhD thesis, Imperial College London, 2010.

\bibitem{kleniati2010partitioning}
P.-M. Kleniati, P.~Parpas, and B.~Rustem.
\newblock Partitioning procedure for polynomial optimization.
\newblock {\em J. Global Optim.}, 48(4):549--567, 2010.

\bibitem{KojimaMuramatsu2009}
M.~Kojima and M.~Muramatsu.
\newblock A note on sparse sos and sdp relaxations for polynomial optimization
  problems over symmetric cones.
\newblock {\em Computational Optimization and Applications}, 42(1):31--41,
  2009.

\bibitem{7285718}
X.~Kuang, L.~F. Zuluaga, B.~Ghaddar, and J.~Naoum-Sawaya.
\newblock Approximating the acopf problem with a hierarchy of socp problems.
\newblock In {\em Power Energy Society General Meeting, 2015 IEEE}, pages 1--5,
  July 2015.

\bibitem{Lasserre1}
J.~Lasserre.
\newblock Global optimization problems with polynomials and the problem of
  moments.
\newblock {\em SIAM Journal on Optimization}, 11:796--817, 2001.

\bibitem{lasserre2015introduction}
J.~B. Lasserre.
\newblock {\em An Introduction to Polynomial and Semi-Algebraic Optimization},
  volume~52.
\newblock Cambridge University Press, 2015.

\bibitem{lavaei2012zero}
J.~Lavaei and S.~Low.
\newblock Zero duality gap in optimal power flow problem.
\newblock {\em IEEE T. Power Syst.}, 27(1):92--107, 2012.

\bibitem{LiuMarecekTakac2015}
J.~{Liu}, A.~C. {Liddell}, J.~{Marecek}, and M.~{Takac}.
\newblock Hybrid methods in solving alternating-current optimal power flows.
\newblock {\em IEEE Transactions on Smart Grid}, 8(6):2988--2998, Nov 2017.

\bibitem{7039413}
A.~Majumdar, A.~A. Ahmadi, and R.~Tedrake.
\newblock Control and verification of high-dimensional systems with dsos and
  sdsos programming.
\newblock In {\em Decision and Control (CDC), 2014 IEEE 53rd Annual Conference
  on}, pages 394--401, Dec 2014.

\bibitem{Marecek2016}
J.~{Marecek}, M.~{Mevissen}, and J.~{Christoffer Villumsen}.
\newblock {MINLP in Transmission Expansion Planning}.
\newblock In {\em 19th Power Systems Computation Conference (PSCC 2016)}, 2016.
\newblock Also, arXiv:1603.04375.

\bibitem{MarecekTakac2015}
J.~Mare{\v c}ek and M.~Tak{\' a}{\v c}.
\newblock A low-rank coordinate-descent algorithm for semidefinite programming
  relaxations of optimal power flow.
\newblock {\em Optimization Methods and Software}, 32(4):849--871, 2017.

\bibitem{Matiyasevich}
Y.~Matiyasevich.
\newblock {\em Hilbert's 10th Problem}.
\newblock The MIT Press, Cambridge, MA, 1993.
\newblock With foreword by Martin Davis and Hilary Putnam.

\bibitem{7038397}
D.~K. Molzahn and I.~A. Hiskens.
\newblock Moment-based relaxation of the optimal power flow problem.
\newblock In {\em Power Systems Computation Conference (PSCC), 2014}, pages
  1--7, Aug 2014.

\bibitem{7232429}
D.~K. Molzahn and I.~A. Hiskens.
\newblock Mixed sdp/socp moment relaxations of the optimal power flow problem.
\newblock In {\em PowerTech, 2015 IEEE Eindhoven}, pages 1--6, June 2015.

\bibitem{MolzahnHiskens2015a}
D.~K. Molzahn and I.~A. Hiskens.
\newblock Sparsity-exploiting moment-based relaxations of the optimal power
  flow problem.
\newblock {\em IEEE Transactions on Power Systems}, 30(6):3168--3180, Nov 2015.

\bibitem{Molzahn2019}
D.~K. Molzahn and I.~A. Hiskens.
\newblock A survey of relaxations and approximations of the power flow
  equations.
\newblock {\em Foundations and Trends® in Electric Energy Systems},
  4(1-2):1--221, 2019.

\bibitem{Powell1969}
M.~J.~D. Powell.
\newblock A method for nonlinear constraints in minimization problems.
\newblock In R.~Fletcher, editor, {\em Optimization}, pages 283--298. Academic
  Press, New York, 1969.

\bibitem{razaviyayn2013}
M.~Razaviyayn, M.~Hong, and Z.-Q. Luo.
\newblock A unified convergence analysis of block successive minimization
  methods for nonsmooth optimization.
\newblock {\em SIAM J. Optimization}, 23(2):1126--1153, 2013.

\bibitem{R76a}
R.~T. Rockafellar.
\newblock Augmented lagrangians and applications of the proximal point
  algorithm in convex programming.
\newblock {\em Math. Oper. Res.}, 1(2):97--116, 1976.

\bibitem{R76b}
R.~T. Rockafellar.
\newblock Monotone operators and the proximal point algorithm.
\newblock {\em SIAM J. Control Optim.}, 14(5):877--898, 1976.

\bibitem{tappenden2013inexact}
R.~Tappenden, P.~Richt{\'a}rik, and J.~Gondzio.
\newblock Inexact coordinate descent: complexity and preconditioning.
\newblock {\em arXiv preprint arXiv:1304.5530}, 2013.

\bibitem{Tawarmalani2005}
M.~Tawarmalani and N.~V. Sahinidis.
\newblock {A polyhedral branch-and-cut approach to global optimization}.
\newblock {\em Mathematical Programming}, 103:225--249, 2005.

\bibitem{vb}
L.~Vandenberghe and S.~Boyd.
\newblock Semidefinite programming.
\newblock {\em SIAM Rev.}, 38(1):49--95, 1996.

\bibitem{Wittek2015}
P.~Wittek.
\newblock Algorithm 950: Ncpol2sdpa\&mdash;sparse semidefinite programming
  relaxations for polynomial optimization problems of noncommuting variables.
\newblock {\em ACM Trans. Math. Softw.}, 41(3):21:1--21:12, June 2015.

\end{thebibliography}
     
\end{document}